\documentclass[12pt,reqno]{article}

\usepackage{amssymb}

\usepackage{fullpage}

\usepackage{amsmath}
\usepackage{amsthm}
\usepackage{amsfonts}
\usepackage{float}

\def\En{{\mathbb N}}

\begin{document}

\theoremstyle{plain}
\newtheorem{theorem}{Theorem}
\newtheorem{corollary}[theorem]{Corollary}
\newtheorem{lemma}[theorem]{Lemma}
\newtheorem{proposition}[theorem]{Proposition}

\theoremstyle{definition}
\newtheorem{definition}[theorem]{Definition}
\newtheorem{example}[theorem]{Example}
\newtheorem{conjecture}[theorem]{Conjecture}
\newtheorem{openproblem}[theorem]{Open Problem}

\theoremstyle{remark}
\newtheorem{remark}[theorem]{Remark}

\def\thre{{3 \over 2}}
\def\mod#1 #2{#1\ ({\rm mod}\ #2)}

\title{Avoiding 3/2-Powers over the Natural Numbers}

\author{Eric Rowland\\
Mathematics Department\\
Tulane University\\
New Orleans, LA 70118\\
USA\\
{\tt erowland@tulane.edu} \\
\and
Jeffrey Shallit\\
School of Computer Science\\
University of Waterloo\\
Waterloo, ON  N2L 3G1\\
Canada\\
{\tt shallit@cs.uwaterloo.ca}
}

\date{January 18, 2011}

\maketitle

\begin{abstract}
In this paper we answer the following question:  what is the 
lexicographically least sequence over the natural numbers that
avoids $\thre$-powers?
\end{abstract}

\section{Introduction}

Ever since the work of Thue more than a hundred years ago, mathematicians
have been interested in avoiding patterns in words.  Thue showed that it 
is possible to create an infinite sequence
over a three-letter alphabet containing
no nonempty squares (that is, no factors of the form $xx$), and
over a two-letter alphabet containing no overlaps (that is, no factors 
of the form $axaxa$, where $a$ is a single letter and $x$ is
a possibly empty word).  

Dejean~\cite{Dejean:1972} instituted the study of fractional powers in 1972.  We say that a
word $x$ is a $p/q$-power, for integers $p > q \geq 1$, if $x$ can be
written in the form $y^e y'$ for some integer $e \geq 1$,
where $y'$ is a prefix of $y$ and $|x|/|y| = p/q$.
For example, the German word {\tt schematische} is a $3/2$-power.

Recently there has been some interest in avoiding patterns over an
infinite alphabet, say $\En$, the natural numbers.  More precisely, we
are interested in finding the {\it lexicographically least} infinite sequence
over $\En$ avoiding $p/q$-powers, by which we mean a sequence
not containing any
factor that is an $\alpha$-power with $\alpha \geq p/q$.
As shown in \cite{Guay-Paquet&Shallit:2009},
this corresponds to running a backtracking
algorithm without actually doing any backtracking; whenever a choice for
the next symbol fails, we increment the choice by $1$ until a valid
choice is found.

Guay-Paquet and the second author \cite{Guay-Paquet&Shallit:2009}
recently described the lexicographically
least words avoiding squares and overlaps over $\En$.  The lexicographically least square-free word over $\En$ is
\[
	01020103010201040102010301020105 \cdots.
\]
Indexing from $0$, the $n$th symbol in this sequence is the exponent of the highest power of $2$ dividing $n+1$.  This sequence is $2$-regular in the sense of Allouche and Shallit \cite{Allouche&Shallit:1992,Allouche&Shallit:2003}.
The lexicographically least overlap-free word can also be described 
compactly, but it is more complicated.

In this paper we are interested in ${\bf w}_{3/2}$, the lexicographically
least word over $\En$ avoiding $\thre$-powers.  As we will see, this word
has a short description, and it is $6$-regular.

Here are the first 100 terms of ${\bf w}_{3/2}$:

\begin{table}[H]
\begin{center}
\begin{tabular}{|r||r|r|r|r|r|r|r|r|r|r|}
\hline
${\bf w}_{3/2}[10i+j]$ & $j = 0$ & $1$ & $2$ & $3$ & $4$ & $5$ & $6$ & $7$ & $8$ & $9$ \\
\hline\hline
$i = 0$ & 0 & 1 & 2 & 0 & 3 & 1 & 0 & 2 & 1 & 3 \\
\hline
$1$ & 0 & 1 & 2 & 0 & 4 & 1 & 0 & 2 & 1 & 4 \\
\hline
$2$ & 0 & 1 & 2 & 0 & 3 & 1 & 0 & 2 & 1 & 5 \\
\hline
$3$ & 0 & 1 & 2 & 0 & 4 & 1 & 0 & 2 & 1 & 3 \\
\hline
$4$ & 0 & 1 & 2 & 0 & 3 & 1 & 0 & 2 & 1 & 4 \\
\hline
$5$ & 0 & 1 & 2 & 0 & 4 & 1 & 0 & 2 & 1 & 5 \\
\hline
$6$ & 0 & 1 & 2 & 0 & 3 & 1 & 0 & 2 & 1 & 3 \\
\hline
$7$ & 0 & 1 & 2 & 0 & 4 & 1 & 0 & 2 & 1 & 4 \\
\hline
$8$ & 0 & 1 & 2 & 0 & 3 & 1 & 0 & 2 & 1 & 6 \\
\hline
$9$ & 0 & 1 & 2 & 0 & 4 & 1 & 0 & 2 & 1 & 3 \\
\hline
\end{tabular}
\end{center}
\end{table}

\section{The word ${\bf w}_{3/2}$}

We define an infinite sequence
${\bf v} = a(0) a(1) a(2) \cdots$
as follows:
\begin{eqnarray}
a(10n) &= & a(10n+3) \; = \; a(10n + 6) \; = \; 0 , \ \ n \geq 0;  \nonumber \\
a(10n+1) &=& a(10n+5) \; = \; a(10n+8) \; = \; 1, \ \ n \geq 0; \nonumber \\
a(10n+2) &=& a(10n+7) \; = \; 2 , \ \ n \geq 0; \nonumber \\
a(10n+4) & = & \begin{cases}
	3, & \text{ if $n \equiv \mod{0} {2}$; } \\
	4, & \text{ if $n \equiv \mod{1} {2}$; } 
	\end{cases} \label{ten} \\
a(10n+9) & = & \begin{cases}
	3, & \text{ if $n \equiv \mod{0} {3}$; } \\
	4, & \text{ if $n \equiv \mod{1} {3}$; } \\
	a(5 (n-2)/3 + 4) + 2, & \text{ if $n \equiv 
		\mod{2} {3}$. }  
	\end{cases} \nonumber
\end{eqnarray}
Thus we have
$${\bf v} = 01203102130120410214012031021501204102130120310214
\cdots \ .$$
Note that ${\bf v}$ has both periodic aspects (in that some linearly-indexed
subsequences are constant) and self-similar aspects (in that
$a(60n+59) = a(10n+9) + 2$).  

Our goal is to show that ${\bf v} = {\bf w}_{3/2}$.  To do so we will prove
that 
(a) ${\bf v}$ is $\thre$-power-free
and
(b) ${\bf v}$ is lexicographically least among all words
avoiding $\thre$-powers.

\section{${\bf v}$ is $\thre$-power-free}\label{power-free}

If $w$ is a $\frac{p}{q}$-power with $\frac{p}{q} \geq \frac{3}{2}$, then $w$ contains a factor of the form $x y x$, where either $|y| = |x|$ or $|y| = |x| - 1$.  We show that neither of these factors appears in $\textbf{v}$.

First, we observe that from the definition, $\bf v$ is
``pseudoperiodic'' with period 10.  More precisely,
$${\bf v} \in
(0120\lbrace 3, 4 \rbrace 1021\lbrace 3, 4, \ldots \rbrace)^\omega ,$$
where by $x^\omega$ we mean the infinite sequence $xxx \cdots$.
From this we immediately see that if $a(n) a(n+1) = a(n+i) a(n+i+1)$ for some $n \geq 0$ and $i \geq 0$ then $i \equiv 0 \pmod {10}$.  It now follows that if $|y| = |x| - 1 \geq 1$ then, since $|xy|$ is odd, $x y x$ is not a factor of $\textbf{v}$.  Similarly, if $|x| = 1$ then $x x$ is not a factor of $\textbf{v}$.

We now show that no word $x y x$, with $|x| = |y| = k \geq 1$, occurs as a factor of $\textbf{v}$.  For each factor $a(n) a(n+1) \cdots a(n + 3k - 1)$ of length $3k$ we exhibit an index $i$, $0 \leq i < k$, such that $a(n + i) \neq a(n + 2k + i)$, implying that this factor is not a $\frac{3}{2}$-power.  The proof is divided up into several cases, depending on the residue class of $k$ modulo $10$.

\bigskip

\noindent{\it Case 1:}  $k \equiv \mod{1, 4, 6, 9} {10}$.
Then from (\ref{ten}) we have $a(n) \not= a(n+2k)$ for all $n$, 
so we can take $i = 0$.

\bigskip

\noindent{\it Case 2:}  $k \equiv \mod{2, 7} {10}$.
If $n \equiv \mod{1, 6} {10}$, then let $i=1$ (which we can do since $k > 1$); otherwise let $i=0$.
One checks that $a(n + i) \neq a(n + 2k + i)$.

\bigskip

\noindent{\it Case 3:}  $k \equiv \mod{3, 8} {10}$.
If $n \equiv \mod{0, 5} {10}$, let $i=1$; otherwise let $i=0$.

\bigskip

\noindent{\it Case 4:}  $k \equiv \mod{5} {10}$.  Then $k \geq 5$.  

If $k = 5$ then
we can choose $i$ such that either $n+i \equiv \mod{4} {10}$
or $n+i \equiv \mod{9} {10}$.  In the former case we have
$a(n+i) = 3$, $a(n+2k+i) = 4$ or vice versa.  In the latter
case we have $a(n+i) \not= a(n+2k+i)$ by (\ref{ten}).

Otherwise $k \geq 15$.  Then we can choose $i$ such that
$n+i \equiv \mod{4} {10}$.  We have $2 k \equiv \mod{10} {20}$, and hence $a(n+i) \not= a(n+2k+i)$ from
(\ref{ten}).

\bigskip

\noindent{\it Case 5:}  $k \equiv \mod{0} {10}$.  
Define
$b(n) = a(10n+9)$ for $n \geq 0$.  Then from (\ref{ten}) we have,
for $n \geq 0$, that
\begin{eqnarray}
b(6n) &=& b(6n+3) \; = \; 3; \nonumber  \\
b(6n+1) &=& b(6n+4) \; = \; 4; \nonumber  \\
b(6n+2) & = & \begin{cases}
	5, & \text{ if $n \equiv \mod{0} {2}$}; \\
	6, & \text{ if $n \equiv \mod{1} {2}$};
	\end{cases} \label{six} \\
b(6n+5) &= & b(n) + 2  \nonumber .
\end{eqnarray}

From this, an easy induction gives
\begin{equation}
b(n) = \begin{cases}
	2t+3, & \text{if the base-$6$ representation of $n$ ends with
		$0 5^t$ or $3 5^t$;} \\
	2t+4, & \text{if the base-$6$ representation of $n$ ends with
		$1 5^t$ or $4 5^t$;} \\
	2t+5, & \text{if the base-$6$ representation of $n$ ends with
		$025^t$, $225^t$ or $425^t$;} \\
	2t+6, & \text{if the base-$6$ representation of $n$ ends with
		$125^t$, $325^t$ or $525^t$.}
	\end{cases}
\label{tee}
\end{equation}

We need to find an index $i$, $0 \leq i < k$,
such that $a(n+i) \not= a(n+2k+i)$.  We choose $i$ such that
$n + i \equiv \mod{9} {10}$.  Since $k = 10r$, there are $r$ possible
choices for $i$.  It follows that we need to show, for each $n \geq 0$
and $r \geq 1$,
that there exists $j$, $0 \leq j < r$ such that
\begin{equation}
b(n+j) \not= b(n+2r+j).
\label{ineq}
\end{equation}

If $r \equiv \mod{1, 2} {3}$, then in fact we can choose $j = 0$ and
use (\ref{six}).

Otherwise $r \equiv \mod{0} {3}$.  This is the most difficult case.
To solve it, we first define two auxiliary sequences, as follows:
\begin{eqnarray*}
c(0) &=& 0; \\
c(6n+1) &=& c(6n+3) \; = \; c(6n+5) \; = \; 2; \\
c(6n+2) &=& c(6n+4) \; = \; 5; \\
c(6n) &=& 6 c(n) + 5;
\end{eqnarray*}
\begin{eqnarray*}
d(0) &=& 0; \\
d(6n+1) &=& d(6n+5) \; = \; 3; \\
d(6n+2) &=& d(6n+3) \; = \; d(6n+4) \; = \; 6; \\
d(6n) &=& 6 d(n) .
\end{eqnarray*}

Write $r = 3s$, and 
let the base-$6$ representation of $s$ end in exactly $t$ zeroes.
An easy induction gives, for $s \geq 1$, that
$$
c(s) = \begin{cases}
	6^{t+1} - 1, & \text{ if the last nonzero digit of $s$ is even;} \\
	3 \cdot 6^t - 1, & \text{ if the last nonzero digit of $s$ is odd;}
	\end{cases}
$$
and
$$
d(s) = \begin{cases}
	6^{t+1}, & \text{ if the last nonzero digit of $s$ is 2, 3, or 4;} \\
	3 \cdot 6^t, & \text{ if the last nonzero digit of $s$ is 1 or 5.}
	\end{cases}
$$

We now claim that for all $s \geq 1$ and $j \geq 0$, we have
\begin{equation}
b(d(s)j + c(s)) \not= b(d(s)j + c(s) + 6s).
\label{bineq}
\end{equation}

Since $c(s) \leq d(s) \leq 3s = r$, (\ref{bineq}) also
provides a solution to (\ref{ineq}).

To verify (\ref{bineq}), we assume that $s$ ends in $t$ 0's.  Then there are
$30$ cases to consider, based on the last nonzero digit $s'$ of $s$
and the last digit $j'$
of $j$ (digit in the base-$6$ representation, of course).
The claim can now be verified by a rather tedious examination of
the 30 cases.  We provide the details for three typical cases:

\bigskip

\noindent $s' = 1$, $j' = 0$:  the base-$6$ expansion of $d(s)j + c(s)$ ends
with $ l 2 5^t$ for some $t \geq 0$ and $l \in \lbrace 0, 3 \rbrace$,
while the base-$6$ expansion of $d(s)j + c(s) + 6s$ ends with
$(l+1) 2 5^t$, so by (\ref{tee}) we see that
$b(d(s)j + c(s))$ and $b(d(s)j+c(s)+6s)$ are
of different parity.

\bigskip

\noindent $s' = 1$, $j' = 3$:
the base-$6$ expansion of $d(s)j + c(s)$ ends
with $1 5^t$ or $4 5^t$ for some $t \geq 0$, while the 
base-$6$ expansion of $d(s)j + c(s) + 6s$ ends with
$2 5^t$ or $5^{t+1}$. From (\ref{tee}) we
get $b(d(s)j + c(s)) = 2t+4$, while
$b(d(s)j+c(s)+6s) \geq 2t+5$.

\bigskip

\noindent $s' = 1$, $j' = 5$:
the base-$6$ expansion of $d(s)j + c(s)$ ends
with $l 5^{t+1}$ for some $l \in \lbrace 2, 5 \rbrace$, while the 
base-$6$ expansion of $d(s)j + c(s) + 6s$ ends with
$((l+1) \bmod 6) 5^{t+1}$.  In either case, from (\ref{tee}) we have
$b(d(s)j + c(s)) \geq 2t+7$, while
$b(d(s)j+c(s)+6s) = 2t+5$.

\bigskip

The proof that $\bf v$ is $\thre$-power free is now complete.

\section{$\bf v$ is the lexicographically least $\thre$-power-free sequence}\label{irreducible}

In this section we complete our characterization of ${\bf w}_{3/2}$
by showing that
a $\thre$-power is formed whenever any symbol of $\bf v$, other than
$0$, is decremented.

The symbol $0$ appears in positions congruent to $0, 3, $ and $6$
(mod $10$), so it suffices to examine other positions.

If the symbol is in a position congruent to $1, 5,$ or $8$ (mod $10$),
then it is a $1$, and it follows $0$, $02$, $03$, or $04$.  Decrementing
this to $0$ then produces the square $00$ or one of the $\thre$-powers $020$, $030$, or $040$.

If the symbol is in a position congruent to $2$ or $7$ (mod $10$),
then it is a $2$, and it follows $01$ or $10$.  Decrementing to
$0$ produces $010$ or $00$, while decrementing to $1$ produces
$11$ or $101$.

If the symbol is in a position congruent to $4$ (mod $10$), then it is
either $3$ or $4$, and it follows $0120$.  Decrementing to $0$ produces
$00$; decrementing to $1$ produces the $5/3$-power $01201$; and decrementing to
$2$ produces $202$.  If the symbol is $4$, there is also the possibility
of decrementing to $3$, and it occurs
at a position $\geq 14$.  If the $4$ is decremented to $3$,
then the immediately preceding
$15$ symbols (including the $3$) form
a $\thre$-power.

The last case is that the symbol is in a position congruent to
$9$ (mod $10$).  Then this position is $10n+9$ for some $n$, and
the symbol is $a(10n+9) = b(n)$.
It therefore follows $1021$.  If it is decremented to $0$, this produces
$10210$.  If it is decremented to $1$, it produces $11$.  If it is
decremented to $2$, it produces $212$.  We now have to handle
the possibility of decrementing to $m$, for some $m$ with
$3 \leq m < b(n)$.

If $n \equiv \mod{0} {3}$, then $b(n) = 3$, so there are no other possibilities
to consider.  If $n \equiv \mod{1} {3}$, then $b(n) = 4$, so we also have
to consider the possibility of decrementing to $3$.  In this case, the
immediately preceding 15 symbols (including the $3$) form a
$\thre$-power.

It remains to consider the case when $n \equiv \mod{2} {3}$.  Here
$b(n)$ is at least $5$.  If we
decrement $b(n)$ to $3$, then the immediately preceding $30$ symbols form a
$\thre$-power, while if we decrement to $4$, then the immediately
preceding $15$ symbols form a $\thre$-power.  

For all other cases, we replace $b(n)$ by $m \geq 5$.  
We claim
this gives the $\thre$-power $xyx$ in $\bf v$, where $x$ is of length
$\ell_m$, in the preceding $3 \ell_m$ symbols, where

$$\ell_m = 
\begin{cases}
30 \cdot 6^{{m \over 2} - 3}, & \text{if $b(n)$ is odd and $m$ is even}; \\
60 \cdot 6^{{{m+1}\over 2} - 3}, & \text{if $b(n)$ is odd and $m$ is odd}; \\
30 \cdot 6^{{m\over 2} - 3}, & \text{if $b(n)$ is even and $m$ is even}; \\
60 \cdot 6^{{{m+1} \over 2} - 3}, & \text{if $b(n)$ is even, $m$ is odd, and
	$m \not= b(n) -1$}; \\
30 \cdot 6^{{{m+1} \over 2} - 3}, & \text{if $b(n)$ is even and $m = b(n) - 1$}.
\end{cases}
$$

To see this, note that from (\ref{ten}) it is enough to show that
\begin{equation}
%	b[n+1-\ell_m/10 \, .. \, n-1] = b[n+1-3\ell_m/10 \, .. \, n-\ell_m/5-1]
	b(n+1-3\ell_m/10) \cdots b(n-\ell_m/5-1) = b(n+1-\ell_m/10) \cdots b(n-1)
\label{b}
\end{equation}
and $b(n-\ell_m/5) = m$.  This can be done by a tedious
examination of each
case in (\ref{tee}).  We give here one representative case,
leaving the rest to the reader.

Suppose the base-$6$ representation of $n$ ends with $3 5^t$ for $t
\geq 1$ and so $b(n) = 2t+3$.

\bigskip

\noindent{\it Case 1:}  $m>6$ is even.  We have $\ell_m/10 = 
3 \cdot 6^{{m \over 2} - 3}$.  The base-$6$ expansion of
$n+1-\ell_m/10$ is $3 5^{t-i-1} 3 0^i$, while the base-$6$
expansion of $n-1$ is $3 5^{t-1} 4$, where $ 1\leq i \leq t-2$ and
$m = 2i+6$.

On the other hand, the base-$6$ expansion
of $n+1-3\ell_m/10$ is
$3 5^{t-i-1} 4 3 0^{i}$, while the base-$6$ expansion of
$n - \ell_m/5 -1$ is
$3 5^{t-i-1} 4 5^{i} 4$.
Using $(\ref{tee})$ we see
that $b$ takes the same values on these intervals.  On the other hand,
$b$ takes the value $m = 2i+6$ at $3 5^{t-i-1} 4 5^{i} 5$.
Thus changing $b(n)$ to $m$ forms a $\thre$-power.

\bigskip

\noindent{\it Case 2:}  $m = 6$.  Then $\ell_m/10 = 3$.  It is now easy to see
that $b(n-5)b(n-4) = b(n-2)b(n-1) = 34$, while $b(n-3) = 6$.
Thus changing $b(n)$ to $m$ forms a $\thre$-power.

\bigskip

\noindent{\it Case 3:}  $m$ is odd.  We have $\ell_m/10 = 6^{{{m+1} \over 2} - 3}$.
The base-$6$ expansion of $n+1-\ell_m/10$ is
$35^{t-i-1} 0^{i+1}$, while the base-$6$ expansion of $n-1$ is
$3 5^{t-1} 4$, where $1 \leq i \leq t-2$ and $m = 2i+5$.

On the other hand, the base-$6$ expansion of
$n+1-3\ell_m/10$ is $3 5^{t-i-2} 3 0^{i+1}$, while the base-$6$
expansion of $n - \ell_m/5 -1$ is
$3 5^{t-i-2} 3 5^i 4$.  Using $(\ref{tee})$ we see that $b$ takes
the same values on these intervals.  On the other hand, $b$
takes the value $m = 2i+5$ at $3 5^{t-i-2} 3 5^{i+1}$.  
Thus changing $b(n)$ to $m$ forms a $\thre$-power.

\section{Morphism description}

The sequence ${\bf w}_{3/2}$ can also be generated as follows.  Consider
the morphisms $\varphi$ and $\tau$ defined by
\begin{eqnarray*}
\varphi(n) &=& 3 \overline{3} 4 \overline{4} 3 (\overline{n+2}) \\
\varphi(\overline{n}) &=& 4 \overline{3} 3 \overline{4} 4 (\overline{n+2}) \\
\tau(n) &=& 0 1 2 0 n \\
\tau(\overline{n}) &=& 1 0 2 1 n.
\end{eqnarray*}
Equation~\eqref{six} is equivalent to the statement that the sequence $3 \, \overline{b(0)} \, 4 \, \overline{b(1)} \, 3 \, \overline{b(2)} \, 4 \, \overline{b(3)} \, \cdots$ is a fixed point of $\varphi$; hence this sequence is $\varphi^\omega(3)$.
It follows that ${\bf w}_{3/2} = \tau(\varphi^\omega(3))$.

\section{Avoiding only $\thre$-powers}

In this section we consider a variant of ${\bf w}_{3/2}$.
Whereas ${\bf w}_{3/2}$ is the lexicographically least word over $\En$
not containing any $\alpha$-power for $\alpha \geq 3/2$, we let ${\bf
x}_{3/2}$ be the lexicographically least word over $\En$ not containing
any (exact) $\thre$-power.  Here are the first $144$ terms of ${\bf
x}_{3/2}$:

\begin{table}[H]
\begin{center}
\begin{tabular}{|r||r|r|r|r|r|r|r|r|r|r|r|r|}
\hline
${\bf x}_{3/2}[12i+j]$ & $j = 0$ & $1$ & $2$ & $3$ & $4$ & $5$ & $6$ & $7$ & $8$ & $9$ & $10$ & $11$ \\
\hline\hline
$i = 0$ & 0 & 0 & 1 & 1 & 0 & 2 & 1 & 0 & 0 & 1 & 1 & 2 \\
\hline
1 & 0 & 0 & 1 & 1 & 0 & 3 & 1 & 0 & 0 & 1 & 1 & 3 \\
\hline
2 & 0 & 0 & 1 & 1 & 0 & 2 & 1 & 0 & 0 & 1 & 1 & 4 \\
\hline
3 & 0 & 0 & 1 & 1 & 0 & 3 & 1 & 0 & 0 & 1 & 1 & 2 \\
\hline
4 & 0 & 0 & 1 & 1 & 0 & 2 & 1 & 0 & 0 & 1 & 1 & 3 \\
\hline
5 & 0 & 0 & 1 & 1 & 0 & 3 & 1 & 0 & 0 & 1 & 1 & 4 \\
\hline
6 & 0 & 0 & 1 & 1 & 0 & 2 & 1 & 0 & 0 & 1 & 1 & 2 \\
\hline
7 & 0 & 0 & 1 & 1 & 0 & 3 & 1 & 0 & 0 & 1 & 1 & 3 \\
\hline
8 & 0 & 0 & 1 & 1 & 0 & 2 & 1 & 0 & 0 & 1 & 1 & 5 \\
\hline
9 & 0 & 0 & 1 & 1 & 0 & 3 & 1 & 0 & 0 & 1 & 1 & 2 \\
\hline
10 & 0 & 0 & 1 & 1 & 0 & 2 & 1 & 0 & 0 & 1 & 1 & 3 \\
\hline
11 & 0 & 0 & 1 & 1 & 0 & 3 & 1 & 0 & 0 & 1 & 1 & 5 \\
\hline
\end{tabular}
\end{center}
\end{table}

Unlike ${\bf w}_{3/2}$, the word ${\bf x}_{3/2}$ contains squares, for example.
However, the underlying structures of ${\bf w}_{3/2}$ and ${\bf x}_{3/2}$ are the same.
Namely, let ${\bf y} = f(0) f(1) f(2) \cdots$, where
\begin{eqnarray*}
f(12n) &= & f(12n+1) \; = \; f(12n+4) \; = \; f(12n+7) \; = \; f(12n+8)\phantom{0} \; = \; 0 , \ \ n \geq 0; \\
f(12n+2) &=& f(10n+3) \; = \; f(12n+6) \; = \; f(12n+9) \; = \; f(12n+10) \; = \; 1, \ \ n \geq 0; \\
f(12n+5) & = & \begin{cases}
	2, & \text{ if $n \equiv \mod{0} {2}$; } \\
	3, & \text{ if $n \equiv \mod{1} {2}$; } 
	\end{cases} \\
f(12n+11) & = & \begin{cases}
	2, & \text{ if $n \equiv \mod{0} {3}$; } \\
	3, & \text{ if $n \equiv \mod{1} {3}$; } \\
	f(2 n + 1) + 2, & \text{ if $n \equiv \mod{2} {3}$. }  
	\end{cases}
\end{eqnarray*}
Note that $f(12 n + 11) + 1 = b(n) = a(10 n + 9)$ for $n \geq 0$.
We show that ${\bf y} = {\bf x}_{3/2}$.
In particular, it follows that ${\bf x}_{3/2}$ is $6$-regular.

First we show that $\bf y$ is $\thre$-power-free.
As in Section~\ref{power-free}, for each $n \geq 0$ and $k \geq 0$ we exhibit an index $i$, $0 \leq i < k$, such that $f(n + i) \neq f(n + 2k + i)$.

\bigskip

\noindent{\it Case 1:}  $k \equiv \mod{1, 5, 7, 11} {12}$.
We have $f(n) \neq f(n+2k)$ for all $n$, so let $i = 0$.

\bigskip

\noindent{\it Case 2:}  $k \equiv \mod{2, 4, 8, 10} {12}$.
If $n$ is even, let $i=1$; if $n$ is odd, let $i=0$.

\bigskip

\noindent{\it Case 3:}  $k \equiv \mod{3, 9} {12}$.
If $n$ is even, let $i=0$; if $n$ is odd, let $i=1$.

\bigskip

\noindent{\it Case 4:}  $k \equiv \mod{6} {12}$.
If $k = 6$ then we can choose $i$ such that either $n+i \equiv \mod{5} {12}$ or $n+i \equiv \mod{11} {12}$; in either case we have $f(n+i) \neq f(n+2k+i)$.
Otherwise $k \geq 18$, and we can choose $i$ such that $n+i \equiv \mod{5} {12}$, and hence $f(n+i) \neq f(n+2k+i)$.

\bigskip

\noindent{\it Case 5:}  $k \equiv \mod{0} {12}$.
Since $f(12 n + 11) + 1 = b(n)$, this case follows immediately from Case~5 of Section~\ref{power-free}.

\bigskip

Therefore $\bf y$ is $\thre$-power-free.
Now we show that decrementing any nonzero symbol of $\bf y$ introduces a $\thre$-power.

The symbol in each position congruent to $0$, $1$, $4$, $7$, or $8$ modulo $12$ is $0$, which cannot be decremented.

The symbol in each position congruent to $2$, $3$, $6$, $9$, or $10$ modulo $12$ is $1$, and it follows $00$, $01$, $02$, or $03$.
In each case, decrementing the $1$ to $0$ produces a $\thre$-power.

The symbol in each position congruent to $5$ modulo $12$ is either $2$ or $3$, and it follows $00110$.
Decrementing to $0$ produces $001100$, and decrementing to $1$ produces $101$.
If the symbol is $3$, then decrementing to $2$ produces a $\thre$-power in the preceding $18$ symbols (including the $2$).

The symbol in each position congruent to $11$ modulo $12$ is at least $2$.
Decrementing to $0$ produces $100110$.
Decrementing to $1$ produces $111$.
If $n \equiv \mod{0} {3}$, then $f(12n+11) = 2$, so there are no other possibilities to consider.
If $n \equiv \mod{1} {3}$, then $f(12n+11) = 3$, and decrementing to $2$ produces a $\thre$-power in the preceding $18$ symbols.

If $n \equiv \mod{2} {3}$, then $f(12n+11) \geq 4$.
Decrementing to $2$ or $3$ produces a $\thre$-power in the preceding $36$ or $18$ symbols, respectively.
It follows from \eqref{b} that if $f(12n + 11) > m \geq 4$ then decrementing $f(12n + 11) = b(n) - 1$ to $m$ produces a $\thre$-power the preceding $\frac{18}{5} \ell_{m+1}$ symbols.

We have shown that $\bf y$ is the lexicographically least $\thre$-power-free word over $\En$; hence ${\bf x}_{3/2} = {\bf y}$.
The word ${\bf x}_{3/2}$ is generated by the morphism $\varphi$ of the previous section as ${\bf x}_{3/2} = \upsilon(\varphi^\omega(3))$, where $\upsilon$ is the morphism defined by
\begin{eqnarray*}
	\upsilon(n) &=& 0 0 1 1 0 (n-1) \\
	\upsilon(\overline{n}) &=& 1 0 0 1 1 (n-1).
\end{eqnarray*}

As already mentioned, ${\bf x}_{3/2}$ contains squares.
However, the only squares in ${\bf x}_{3/2}$ are $00$ and $11$, as we now show.
Suppose $xx$ is a factor of ${\bf x}_{3/2}$ for some nonempty word $x$.
The length of $x$ cannot be even, because otherwise $xx$ contains a $\thre$-power.
One checks from the definition of $f(n)$ that the length of $x$ cannot be $3$.
Similarly, $f(n) \cdots f(n + 4) \neq f(n + k) \cdots f(n + 4 + k)$ for all $n \geq 0$ and odd $k$, which precludes a factor $xx$ where the length of $x$ is odd and at least $5$.

Hence the only square factors of ${\bf x}_{3/2}$ are $00$ and $11$.
Since $000$ and $111$ are not factors of ${\bf x}_{3/2}$, it follows that ${\bf x}_{3/2}$ is overlap-free.

% Since ${\bf w}_{3/2}$ avoids $\alpha$-powers for each $\alpha \geq
% 3/2$, it follows that ${\bf x}_{3/2}$ avoids sufficiently long
% $\alpha$-powers for each $\alpha \geq 3/2$.

% Conjecturally, the set of $\alpha > 3/2$ such that ${\bf x}_{3/2}$
% contains an $\alpha$-power is $\{\frac{11}{7}, \frac{8}{5},
% \frac{5}{3}, \frac{9}{5}, 2\}$.

\section{Remarks}

A previous paper discussed the structure of the lexicographically
least word over $\En$ avoiding $n$'th powers or overlaps
\cite{Guay-Paquet&Shallit:2009}.  In this paper, we discussed the
structure of such a word avoiding $3/2$-powers.  It remains to gain
a deeper understanding of the lexicographically least word avoiding
$\alpha$-powers for arbitrary $\alpha$.  Computer experiments strongly
suggest that a similar structure exists for $4/3$-powers.  But it is
still not known whether the lexicographically least word avoiding $5/2$-powers
uses only finitely many distinct letters.

\bibliographystyle{plain}

\begin{thebibliography}{9}

\bibitem{Allouche&Shallit:1992}
J.-P. Allouche and J. Shallit.
\newblock The ring of $k$-regular sequences.
\newblock {\it Theoret. Comput. Sci.} {\bf 98} (1992), 163--197.

\bibitem{Allouche&Shallit:2003}
J.-P. Allouche and J. Shallit.
\newblock The ring of $k$-regular sequences, II.
\newblock {\it Theoret. Comput. Sci.} {\bf 307} (2003), 3--29.

\bibitem{Dejean:1972}
F. Dejean.
\newblock Sur un {th\'eor\`eme} de {Thue}.
\newblock {\it J. Combin. Theory Ser. A} {\bf 13} (1972), 90--99.

\bibitem{Guay-Paquet&Shallit:2009}
M. Guay-Paquet and J. Shallit.
\newblock Avoiding squares and overlaps over the natural numbers.
\newblock {\it Discrete Math.} {\bf 309} (2009), 6245--6254.

\end{thebibliography}

\end{document}